\input amstex      
\documentstyle{amsppt}      
      
\magnification=1200      
\pagewidth{6.5truein}      
\pageheight{9.0truein}

\define\C{\Bbb C}      
    
\define\R{\Bbb R}      
      
\define\N{\Bbb{N}}
      
\define\dbar{\bar\partial}      
      
\define\lam{\lambda}

\define\sumprime{\sideset \and^\prime \to \sum}

\redefine\Re{\operatorname{Re}}

\topmatter      
\title      
Compactness of the $\dbar$-Neumann Problem on Convex Domains
\endtitle      
\rightheadtext{Compactness of the $\dbar$-Neumann problem}      
\author      
Siqi Fu and Emil J. Straube
\endauthor      
      
%\affil      
%Texas A \& M University 
%\endaffil      
      
\address
Department of Mathematics, \ \ Texas A \& M University,
College Station, \ \ TX 77843
\endaddress
      
\email       
sfu\@math.tamu.edu \ \ straube\@math.tamu.edu 
\endemail

\thanks
Research supported in part by NSF grant DMS-9500916
\endthanks

%\date      
%November 1997
%\enddate      
      

\abstract 
The $\dbar$-Neumann operator on $(0, q)$-forms 
($1\le q \le n$) on a bounded
convex domain $\Omega$ in $\C^n$ is compact if and only if 
the boundary of $\Omega$ contains no complex analytic (equivalently: 
affine) variety
of dimension greater than or equal to $q$.      
\endabstract

\endtopmatter

\document      
    
\subhead 1. Introduction \endsubhead    

\smallskip

Let $\Omega$ denote a bounded pseudoconvex domain in $\C^n$.
For $1\le q\le n$,  let $L^2_{(0, q)}(\Omega)$
denote the space of $(0, q)$-forms with square integrable 
coefficients, with the norm $\|\sum^\prime a_{I}d\bar z_I\|^2$
$=\sum^\prime \int_\Omega |a_{I}|^2 dV(z)$, where the prime
indicates the summation over strictly increasing $q$-tuples.
The $\dbar$-Neumann operator $N_q$ is the (bounded) inverse
of the (unbounded) self-adjoint, surjective operator 
$\dbar\dbar^* +\dbar^*\dbar$.  We refer the reader to [FK],
[Ko], [Kr2], and the recent survey [BS3] for background on the 
$\dbar$-Neumann problem.

Compactness of the $\dbar$-Neumann problem is a basic property
with many useful consequences. In the case of domains with smooth
boundary, it implies global regularity of the $\dbar$-Neumann
problem (in the sense of preservation of the $L^2$-Sobolev spaces),
see [KN]. Also, the Fredholm theory for Toeplitz operators
is a direct consequence of the compactness of the $\dbar$-Neumann
problem ([V], [HI]). In a related context, whether or not the
$\dbar$-Neumann problem is compact has ramifications for certain
$C^*$-algebras of operators naturally associated with a domain
in $\C^n$; compare for example [Sa]. 

Catlin [Ca2] proved compactness of the $\dbar$-Neumann operator
on smoothly bounded domains whose boundary satisfies Property (P).
The boundary of a domain $\Omega$ satisfies Property (P) if for
every positive number $M$ there is a plurisubharmonic function
$\lam\in C^\infty(\overline{\Omega})$ with $0\le \lam\le 1$,
such that for all $z\in b\Omega$ and $w\in \C^n$,
$\sum_{\alpha, \beta=1}^n
(\partial^2 \lam(z)/\partial z_\alpha \partial \bar z_\beta)w_\alpha
\bar w_\beta \ge M|w|^2$ ([Ca2]).
Property (P) was studied systematically under the name of $B$-regularity
in [Si] (in the context of arbitrary compact sets in $\C^n$). 
It was recently observed that
Catlin's result remains true when no boundary smoothness at all 
is assumed: the $\dbar$-Neumann problem is compact on a bounded
pseudoconvex domain whose boundary is a $B$-regular set ([St]). 
(Under the additional assumption that the domain is hyperconvex,
compactness of the $\dbar$-Neumann problem had been shown earlier
in [HI].)

In this article, we discuss compactness of the $\dbar$-Neumann
problem on bounded convex domains. We obtain a complete
characterization of compactness by the absence from the
boundary of complex analytic (equivalently: affine)
varieties of appropriate dimensions.
A closely related (in fact, in the context of convex domains,
equivalent) question is that of compact solution operators for
$\dbar$. For $1\le q\le n$, 
consider $\dbar$ as an unbounded
operator from $L^2_{(0, q-1)}(\Omega)$ to $L^2_{(0, q)}\cap
\ker\dbar$.  A bounded linear operator $S_q$ from $L^2_{(0, q)}\cap
\ker\dbar$ to  $L^2_{(0, q-1)}(\Omega)$ is called a solution
operator for $\dbar$ on $(0, q)$-forms if $\dbar S_q u =u$
for all $u\in L^2_{(0, q)}\cap \ker\dbar$. 

The following terminology
will be convenient: an affine variety of dimension $q$
is a (relatively) open
subset of a complex affine subspace of $\C^n$ of dimension $q$.
 
\proclaim{Theorem 1.1} Let $\Omega$ be a bounded convex domain
in $\C^n$. Let $1\le q\le n$.  The following are
equivalent:
\itemitem{(1)} There exists a compact solution operator for $\dbar$
on $(0, q)$-forms.
\itemitem{(2)} The boundary of $\Omega$ does not contain
any affine variety of dimension greater than or equal to $q$.
\itemitem{(3)} The boundary of $\Omega$ does not contain
any analytic variety of dimension greater than or equal to $q$.
\itemitem{(4)} The $\dbar$-Neumann operator $N_q$ 
is compact.
\endproclaim

The implication (4) $\Rightarrow$ (1) holds in general (and
is well known): compactness of the $\dbar$-Neumann operator
implies compactness of the canonical solution operator. In fact,
the formula $N_q=(\dbar^*N_q)^*(\dbar^* N_q) + (\dbar^*N_{q+1})
(\dbar^*N_{q+1})^*$ (see [R], [FK]) shows that $N_q$ is compact
if and only if the canonical solution operators $\dbar^*N_q$ and $\dbar^*
N_{q+1}$ are compact. 
Note that in statement (1), it is
the same to say that the canonical solution operator is compact,
since compactness is preserved by projection onto the orthogonal
complement of the kernel of $\dbar$. 

Henkin and Iordan [HI]
recently showed that (on a bounded convex domain) the 
$\dbar$-Neumann operators are compact
for $1\le q \le n$ if there are no one-dimensional
analytic varieties contained in the boundary of the domain
(see [Si] for the smooth case).  
It has  been known  for some time that analytic discs in 
the boundary of a smooth pseudoconvex domain in $\C^2$
obstruct compactness of the $\dbar$-Neumann operator. 
Specific examples for failure of compactness 
are given in [Li] and [Kr2]. These examples are pseudoconvex
Reinhardt domains in $\C^n$.  Theorem 5 in [Sa] implies that on a
pseudoconvex Reinhardt domain, compactness of 
the canonical solution operator on $(0, 1)$-forms 
is incompatible with analytic discs on 
the boundary.  To what extent analytic varieties in the
boundary obstruct compactness of the $\dbar$-Neumann problem on 
``general'' domains seems to be open.   On the other hand, it is 
known that obstructions to compactness of the $\dbar$-Neumann problem
can be more subtle than analytic varieties in the boundary
(so that the above characterization is false without some assumption
on the domain): Matheos [Mt] has recently shown that there exist smooth
bounded pseudoconvex Hartogs domains in $\C^2$ without discs in their
boundaries, whose $\dbar$-Neumann operators are nonetheless not compact.

We remark that we make no explicit assumption on smoothness of
the boundary. However, convexity implies that the boundary is
Lipschitz (see e.g. [Mz], \S 1.1.8).  

The remainder of the  paper is organized as follows.  
In Section 2, we briefly discuss the equivalence of conditions 
(2) and (3) in Theorem 1.1.  This may be part of the folklore,
but we include the (simple) argument for completeness. Section 3
contains the proof that (2) implies (4). The proof of Theorem 1.1
is completed in Section 4, where we show that (1) implies (2).
We conclude the paper with some additional remarks in Section 5.

\bigskip

\noindent{ Acknowledgement:}\ The authors thank Harold Boas for
stimulating conversations.

\bigskip

\subhead 2. Varieties in the Boundary of a Convex Domain \endsubhead

\smallskip

In this section we show that (2) implies (3) in Theorem 1.1 ((3)$\Rightarrow$
(2) is trivial).  This may be viewed as the simplest manifestation of the 
general principle that on the boundaries of convex domains, questions of 
orders of contact with analytic varieties are determined by the orders
of contact with affine subspaces ([Mc], [BS2], [Y]).
Related observations may be found in [N], [Ch].
Complex manifolds in general (but smooth) pseudoconvex boundaries are
studied in [BF].

We first observe the following: if $V$ is a $q$-dimensional variety in 
$\C^n$, then its convex hull $\widehat V$
contains an affine variety of dimension $q$.
This is clear if $n=1$.  For general $n$, the observation follows by induction
on the dimensions as follows. If $\widehat V$ has non-empty interior (in
$\C^n$), we are done. If the interior of $\widehat V$ is empty, $V$ is
contained in a real hyperplane (since then there are no $2n$ line
segments with  end points in $V$ which are linearly independent over
$\R$).  After a suitable change of coordinates, this hyperplane
is $\{ x_n=0\}$.  By the open mapping property of non-constant
holomorphic functions, applied to the restriction of the function 
$z_n$ to (the regular part of) $V$, $V$  is contained in the complex
hyperplane $\{z_n=0 \}$.  This completes the induction.

To prove that (2) implies (3) in Theorem 1.1, assume now that $b\Omega$
contains a $q$-dimensional analytic variety $V$.  Let $p_0$ be a regular
point of $V$ so that near $p_0$, $V$ is a $q$-dimensional complex manifold,
and assume without loss of generality that a supporting hyperplane for 
$\Omega$ at $p_0$ is given by $\{x_n=0\}$.  The argument in the previous
paragraph shows that if $V_{p_0}$ is the intersection of $V$ with a small
neighborhood of $p_0$, then $V_{p_0}\subseteq \{z_n=0\}$.   
Consequently, the convex hull of $V_{p_0}$ is likewise contained both in 
$\overline{\Omega}$ and in $\{z_n=0\}\subseteq \{x_n=0\}$, hence in
$b\Omega$.  But by the above observation, this convex hull contains
a $q$-dimensional affine variety.

\bigskip

\subhead 3. Sufficient Conditions for  Compactness of $N_q$\endsubhead
\smallskip

\proclaim{Proposition 3.1} Let $\Omega$ be a bounded pseudoconvex
domain in $\C^n$.
Assume that for every positive number $M$, there exists a neighborhood 
$U$ of $b\Omega$ and a $C^2$-smooth function $\lam$ on $U$, $0\le \lam \le 1$,
such that for all $z\in U$, the sum of any $q$ (equivalently: the smallest
$q$) eigenvalues of the Hermitian form $(\partial^2\lam(z)/\partial z_\alpha
\partial \bar{z}_\beta)_{\alpha, \beta=1}^n$ is at least $M$.
%$$
%\sumprime_{|I|=q-1} 
%\sum^n_{\alpha, \beta=1}
%\frac{\partial^2\lambda(z)}{\partial z_\alpha\partial \bar{z}_\beta}
%f_{\alpha I}\overline{f_{\beta I}} \ge M\|f\|^2
%\tag 3.1
%$$
%for all $z\in U$ and $f\in \Lambda^{(0, q)}_z$ , 
%with the prime indicating the sum over all strictly increasing 
%$(q-1)$-tuples $I$.  
Then the $\dbar$-Neumann operator $N_q$ on $\Omega$ is compact.
\endproclaim

Note that for $q=1$, the above condition on the Hessian
of $\lam$ reduces to the condition that appears
in Property (P).  For $q>1$, this condition 
does not imply that $\lam$ is {\it pluri}subharmonic
(while  it still implies that $\lam$ is subharmonic).

To prove Proposition 3.1, we first note (following
[H], p.~137) that the condition on $\lam$ implies 
$$
\sumprime_{I} 
\sum^n_{\alpha, \beta=1}
\frac{\partial^2\lambda(z)}{\partial z_\alpha\partial \bar{z}_\beta}
f_{\alpha I}\overline{f_{\beta I}} \ge M\|f\|^2
\tag 3.1
$$
for all $z\in U$ and $f\in \Lambda^{(0, q)}_z$, where
$\Lambda^{(0, q)}_z$ denotes the space of $(0, q)$-forms
at $z$, and the prime indicates summation over increasing 
$(q-1)$-tuples $I$. (3.1) can be seen by using a frame (at $z$)
where the Hessian of $\lam$ is diagonalized. The proof of Proposition
3.1 is now based on 
[Ca3], Theorem 2.1 (see also [BS3], Section 2,
for a somewhat different approach to this type of estimate) and
[Ca2], proof of Theorem 1.  The fact that no boundary smoothness
is assumed necessitates working on smooth subdomains and using
a regularization procedure for the forms involved that was
introduced in [St]. The details of this argument are carried out in [St],
proof of Corollary 3, to which we refer the reader.

\medskip

In order to show that (2)$\Rightarrow$(4) in Theorem 1.1, it now
suffices to show that if the boundary of a bounded convex domain
contains no affine varieties of dimension $q$ or greater, then 
the assumption in Proposition 3.1 is satisfied.  This can be done
by suitably generalizing the arguments in [Si], Proposition 2.4
that cover the case $q=1$: the Choquet theory has to be done
for a cone of functions that reflects the condition on the Hessian used in Proposition 3.1 (rather than
for the cone of plurisubharmonic functions).  
We also need a substitute for the fact that boundary points
of a convex domain $\Omega$  are peak points for the algebra of
functions holomorphic on $\Omega$ and continuous on $\overline{\Omega}$
when there are no varieties of positive dimension in the boundary
of $\Omega$.  We now develop the necessary ideas.

For a compact set $X$ in $\C^n$, denote by $C(X)$ the
usual algebra of continuous functions on $X$, with 
sup-norm. Denote by $H(X)$ the closure in $C(X)$
of the functions holomorphic in a neighborhood of $X$.
A closed subset $E$ of $X$ is said
to be a peak set if there exists $f\in H(X)$ such that $f(z)=1$
on  $E$ and $|f(z)|<1$ on $X\setminus E$. The function
$f$ is called a (weak) peak function on $E$.

\proclaim{Proposition 3.2} Let $X$ be a compact convex subset of $\C^n$,
let $z_0\in X$, and let $1\le q\le n$.  Then there exists a complex affine
subspace $L$
of dimension $\le q-1$ through $z_0$ such that $X\cap L$ is a peak 
set if and only if $X$ contains no affine variety of dimension $\ge q$
through $z_0$.
\endproclaim

\demo{Proof} The $\Rightarrow$ direction follows easily from
the maximum modulus principle.  To prove the reverse direction,
we need Glicksberg's peak set theorem (cf. [G2], pp.~58) which
says that $E$ is a peak set if and only if $\nu_E\in H(X)^\perp$ for
all finite regular Borel measures $\nu\in H(X)^\perp$, where $\nu_E$ is the
restriction of $\nu$ to $E$.

We argue by induction.  The case $n=1$ is clear (but see the case $q=n$
below).  Assume the conclusion for $n-1$, ($n\ge 2$). 
We need to establish it for $n$.  Since there is no affine 
variety of dimension $\ge q$ through $z_0$ and contained in $X$,
$z_0$ is a boundary point.  Without loss of
generality, assume that $z_0$ is the origin and 
$X\subseteq \{\Re z_n\ge 0\}$.
Set $g(z)=\exp(-\sqrt{z_n})$ (where the square root is the principal branch).
Let $J=X\cap \{z_n=0\}$.  If $q=n$, let $L=\{z_n=0\}$. Then $J$ is a peak
set: the function $g(z)$ is a peak function.  Now assume that
$1\le q\le n-1$. Then $J$ is a compact convex subset of $\C^{n-1}$, and 
by the induction assumption, there is a complex affine
subspace  $L\subseteq \{z_n =0\}\cong \C^{n-1}$ 
of dimension $\le q-1$ such that $J\cap L$ is a peak set for $H(J)$.  
We now show that $X\cap L=J\cap L$ is a peak set for $H(X)$.
Let $\nu$ be a finite regular Borel measure on $X$, $\nu \in H(X)^\perp$.  
Thus we have 
for any holomorphic polynomial $f$ and positive integer $m$ that
$\int_X f\cdot g^m \, d\nu=0$ 
(note that $g$, although not itself analytic in a neighborhood
of $X$, is in $H(X)$).
Letting $m\to \infty$, we obtain that $\int_J f\, d\nu=0$.  
The convex set $J$ is polynomially convex, so 
the holomorphic polynomials are dense in $H(J)$  by the Oka-Weil
approximation theorem. Consequently,  
$\nu_J\in H(J)^\perp$ and hence $\nu_{J\cap L}\in H(J)^\perp$,
by Glicksberg's theorem.  Using Glicksberg's theorem in the
other direction, we conclude that $J\cap L$ is a peak set for
$H(X)$. This completes the induction and the proof of Proposition 3.2.
\enddemo

For an open set $U\subset \C^n$, denote by $P_q(U)$
the set of continuous functions $\lam$ on $U$ such that
for any $z\in U$ and orthonormal set of vectors 
$\{t_1, \cdots, t_q\}$ 
in $\C^n$, the function 
$$
\zeta=(\zeta_1,\cdots, \zeta_q)\in\C^q
\mapsto \lam (z + \zeta_1 t_1
+\ldots +\zeta_q t_q)
$$ 
is subharmonic on $\{\zeta\in \C^q;\ 
z+\zeta_1 t_1 +\ldots + \zeta_q t_q\in U \}$.  
That is, $P_q(U)$ consists of the continuous functions
on $U$ that are subharmonic on each $q$-dimensional complex affine
subspace.  In particular, $P_1(U)$ is the set of all continuous
plurisubharmonic functions and $P_n(U)$ is the set of
all continuous subharmonic functions. 
%A $C^2$-smooth
%function $\lam$ is in $P_q(U)$ if and only if 
%the sum of any $q$ eigenvalues of the Hermitian matrix
%$(\partial^2 \lam (z)/\partial z_\alpha \dbar z_\beta)$
%is non-negative for any $z\in U$. See the remarks after Proposition~3.1
%for other equivalent conditions. 
$P_q(U)$ is a convex cone in $C(U)$ that is closed under taking the pointwise
maximum of finitely many of its elements.  Note that each
function in $P_q(U)$ is a locally uniform limit of $C^\infty$-smooth
elements in $P_q$ of slightly smaller open sets: this follows
from the usual mollifier argument.  Finally, it is
not hard to check that $-\sum_{j=1}^{q-1}|z_j|^2
+(q-1)\sum_{j=q}^n |z_j|^2\in P_q(\C^n)$.

We now return to a bounded pseudoconvex domain $\Omega$.
We denote by $P_q(b\Omega)$ the closure in $C(b\Omega)$
of functions that are in $P_q$ in a neighborhood of $b\Omega$.
A probability measure $\mu$ on $b\Omega$ is said to be
a $P_q$-measure for $z\in b\Omega$ if
$$
\lambda(z)\le \int_{b\Omega} 
\lambda\, d\mu, \qquad \lambda\in P_q(b\Omega). \tag 3.2
$$
We refer the reader to [G1], Chapter 1 for a treatment of these
measures in an abstract context, and for the elements of 
Choquet theory.  In particular, $P_q(b\Omega)$ satisfies
the properties (1.1)-(1.3) in [G1].

Let $\Omega$ be a bounded convex domain, and $z_0$ a boundary
point through which there is no affine variety of dimension $q$
or higher that is contained in $b\Omega$.  We claim that 
the only $P_q$-measure for $z_0$ is the point mass at $z_0$.
Note that there is also no affine variety of dimension
$\ge q$ through $z_0$ that is contained in $\overline{\Omega}$
(this is a special case of the argument at the end of Section 2).
By Proposition 3.2, 
there is a complex affine subspace $L$
of dimension $\le q-1$ such that $b\Omega\cap L$
is a peak set for $H(\overline{\Omega})$.  
Let $f$ be the corresponding (weak)
peak function. Because $f\in H(\overline{\Omega})$,
$|f|\in P_q(b\Omega)$. (3.2) now shows that
any $P_q$-measure $\nu$ for $z_0$ is supported on $b\Omega\cap L$.
In suitable coordinates, we may assume that
$z_0$ is the origin and $L\subseteq \{z_q=\ldots =z_n=0 \}$.  In (3.2),
take now $\lambda = -\sum_{j=1}^{q-1} |z_j|^2
+(q-1)\sum_{j=q}^n |z_j|^2$.  We already know that the support
of $\nu$ is contained in $L$, where $\lam$ reduces 
to $-\sum_{j=1}^{q-1} |z_j|^2$.  We thus obtain from 
(3.2) that the support of $\nu$ consists of the point $z_0$. 

We next invoke Edwards' theorem ([G1], Theorem 1.2): for every continuous function $u$ on $b\Omega$ and $z\in b\Omega$, $\inf\{\int_{b\Omega}u\,d\mu\,;\enspace \mu \text{ is a $P_{q}$-measure for $z$}\} = \sup\{\lambda(z);\enspace \lambda\in P_q(b\Omega), 
\lambda \le u \text{ on }
 \ b\Omega \}$. Because all
$P_q$-measures have point support, the theorem gives
$$
u(z)=\sup\{\lambda(z);\enspace \lambda\in P_q(b\Omega), 
\lambda \le u \text{ on }
 \ b\Omega \} \tag 3.3
$$
for every function $u\in C(b\Omega)$.
For $M>0$, let $u_M(z)= -M|z|^2$.  It follows from (3.3) and
a compactness argument similar to the proof of Dini's theorem
that $u_M$ can be approximated uniformly on $b\Omega$
by functions in $P_q(b\Omega)$, hence by functions
that are smooth and in $P_q$ in a neighborhood 
of $b\Omega$. In particular, there exists a neighborhood
$U$ of $b\Omega$  and a function $\lam\in P_q(U)\cap C^2(U)$
such that $0\le \lam + M|z|^2 \le 1$ on $U$ (after shrinking
$U$ if necessary).  The sum of the smallest $q$ eigenvalues
of the Hessian of $\lam + M|z|^2$ is at least $qM\ge M$ (because
$\lam\in P_q(U)$, the sum of the $q$ smallest eigenvalues
of the Hessian of $\lam$ is non-negative). Therefore $\Omega$
satisfies the assumptions in Proposition 3.1, and the proof 
that (2) implies (4) in Theorem 1.1 is complete.  

\bigskip

\subhead  4. Necessary Conditions for the Compactness of  
Solution Operators to $\dbar$  \endsubhead
\smallskip

In this section, we prove the implication (1)$\Rightarrow$(2) of 
Theorem 1.1.  One of
the main tools in the proof is the Ohsawa-Takegoshi extension theorem 
([OT], [O2]).
We also use an idea that comes from [Ca1] (see also [DP]). 

We first prove an auxiliary lemma. Denote by
$K_\Omega(z, w)$ the Bergman kernel function
of a domain~$\Omega$.

\proclaim{Lemma 4.1} Let $\Omega$ be a bounded convex domain in $\C^n$.
\itemitem{(1)} For any $p_0\in b\Omega$ and $p_1\in \Omega$, there
exist constants $C>0$ and $\delta_0>0$ such that 
$$
K_\Omega(p_\delta, p_\delta)\ge C K_\Omega(p_{2\delta}, 
p_{2\delta}) \tag 4.1
$$
for any $\delta \in (0, \delta_0)$, where $p_\delta=p_0 +
\delta(p_1-p_0)/\|p_1-p_0\|$.

%\itemitem{(2)} There exists a constant $C>0$, depending only on
%the diameter of $\Omega$, such that 
%$$
%K_\Omega(z, z)\ge C\frac{1}{d^2_\Omega(z)} \tag 3.5 $$
%for all $z\in \Omega$, where $d_\Omega(z)$ is the Euclidean 
%distance from $z$ to $b\Omega$.

\itemitem{(2)} For any sequence $\{p_j\}\in \Omega$ converging to 
$p_0\in b\Omega$,
$$
\lim_{j\to \infty} \frac{K_\Omega(z, p_j)}{\sqrt{K_\Omega(p_j, p_j)}} =0, 
\tag 4.2  
$$ 
locally uniformly on $\Omega$.
\endproclaim

\demo{Proof} (1) Let $U$ be a ball with center $p_0$  and radius $r$ 
the minimum of  
$\text{dist}(p_1, b\Omega)$ and $\|p_1-p_0\|/2$.  
Let $\vec{n}=(p_1-p_0)/\|p_1-p_0\|$. It is easy to see
from the convexity of $\Omega$ that $T_\delta(\Omega\cap U)\subseteq
\Omega$ for $0<\delta < \|p_1-p_0\|/2$, where
$T_\delta (z)=z+\delta\vec{n}$. Let $\delta_0 =r/2$.

Then  for $0\le \delta\le \delta_0$,
$$
K_{\Omega}(p_\delta, p_\delta)\ge C K_{\Omega\cap U}(p_\delta, p_\delta) 
=C K_{T_\delta(\Omega\cap U)}(p_{2\delta}, p_{2\delta})
\ge C K_\Omega(p_{2\delta}, p_{2\delta}),
$$
where the first inequality follows by localization
of the kernel ([JP], Theorem 6.3.5), and the last 
inequality holds because $T_\delta(\Omega\cap U)\subseteq
\Omega$. 

\smallskip

(2) This part of the lemma is implicit in work of Pflug (see
[JP, \S 7.6]) and Ohsawa [O1] on the completeness
of the Bergman metric. We recall  the proof for the reader's
convenience.  
Without loss of generality, assume that $\Omega$ contains the origin.
It suffices to establish pointwise convergence: Vitali's
theorem  (note that $K_\Omega(\cdot, p_j)/\sqrt{K_\Omega(p_j, p_j)}$ has 
norm 1) then implies that the convergence is locally uniform.
For $z\in\Omega$, let $f(w)=K_\Omega(z, w)$.  
Then $\|f(w)-f(tw)\|_\Omega \to 0$ 
as $t\to 1^-$. Now fix $t$, $0<t<1$. Then
$$\align
\frac{|f(p_j)|}{\sqrt{K_\Omega(p_j, p_j)}} &\le 
\frac{|f(p_j)-f(tp_j)|}{\sqrt{K_\Omega(p_j, p_j)}}
+\frac{|f(tp_j)|}{\sqrt{K_\Omega(p_j, p_j)}} \\      
&\le \|f(w)-f(tw)\|_\Omega + \frac{|K_\Omega(z, tp_j)|}
{\sqrt{K_\Omega(p_j, p_j)}} .
\endalign
$$
The domain $\Omega$ is convex and 
so satisfies an outer cone condition.  Therefore,
$K_\Omega(p_j, p_j)\to \infty$
as $j\to \infty$ (see e.g. [JP], Theorem 6.1.17).
Thus, letting first $j\to \infty$, then $t\to 1^-$, we 
obtain part (2) of Lemma 4.1. 
\enddemo

We are now in a position to prove the implication (1)$\Rightarrow$(2)
in Theorem 1.1.

\demo{Proof of (1)$\Rightarrow$(2) } 
Arguing indirectly, we assume that there exists
a compact solution operator $S_q$ on $(0, q)$-forms and
$b\Omega$ contains an affine variety of dimension
$q$. (Thus $q\le n-1$.)  After an affine transformation, we
may assume that $\{(z', 0)\in \C^n; \ \ |z'|<2\}\subseteq b\Omega$,
where $z'=(z_1, \ldots, z_q)$.  Let $z''=(z_{q+1}, \ldots, z_n)$.

Let $\Omega_1=\{z''\in \C^{n-q}, \ (0, z'')\in \Omega\}$. 
It follows from the convexity
of $\Omega$ that $\Omega_1$ is a (non-empty) convex domain 
in $\C^{n-q}(z'')$. Let $\Omega_2=\{z''\in \C^{n-q}; \ \ 
2z''\in \Omega_1 \}$.  
Then  $\{z'\in \C^q; \ \ |z'|<1\}\times \Omega_2 \subseteq
\Omega$: every point in this set is the midpoint of a line segment
joining a point in $\{|z'|<2\}\times \{0\}$ to a point in $\{0\}\times
\Omega_1$.

Let $p_0$ be a point in $\Omega_2$ and let $p_j=p_0/j, j\in \N$.  Let
$$
f_j(z'')=\frac{K_{\Omega_1}(z'', p_j)}{\sqrt{K_{\Omega_1}(p_j, p_j)}}.
$$
Then $\|f_j\|_{\Omega_1}=1$. We have  
$$
\aligned
\|f_j(z'')\|^2_{\Omega_2} &=\frac{\|K_{\Omega_1}(\cdot, p_j)\|^2_{\Omega_2}}
{K_{\Omega_1}(p_j, p_j)} 
\ge \frac{K_{\Omega_1}(p_j, p_j)}
{K_{\Omega_2}(p_j, p_j)} \\
&=2^{-2(n-q)}\frac{K_{\Omega_1}(p_j, p_j)}{K_{\Omega_1}(2p_j, 2p_j)}
\ge C ,
\endaligned 
$$
for $j$ large enough. The first inequality follows because
$K_{\Omega_1}(p_j, p_j)\le (K_{\Omega_2}(p_j, p_j))^{1/2}$ $
\|K_{\Omega_1}(\cdot, p_j)\|_{\Omega_2}$ ( obtained by applying the reproducing property of $K_{\Omega_2}(p_j,\cdot)$ to the function $K_{\Omega_1}(\cdot,p_j)$). 
The last equality follows from the transformation formula of the 
Bergman kernel. The last inequality follows from (4.1). 
On the other hand, by (4.2), $f_j\to 0$ locally uniformly
on $\Omega_1$. Consequently, no subsequence of $\{f_j\}$
can converge in $L^2(\Omega_2)$.

By the Ohsawa-Takegoshi extension theorem [OT] (see also [O2]),
there exist $L^2$-holomorphic functions $F_j(z', z'')$ on $\Omega$
such that $F_j(0, z'')=f_j(z'')$ and $\|F_j\|_\Omega \le C$. 
We now use an idea from [Ca1] (compare also [DP]). Let
$\alpha_j=F_j(z', z'') d\bar z_1\wedge\cdots\wedge d\bar z_q$.
Then $\dbar\alpha_j=0$, $\|\alpha_j\|_{L^2_{(0, q)}(\Omega)} \le C$.  
Let $g_j=S_q \alpha_j$. Denote by $\widehat g_j$ the form obtained
from $g_j$ by discarding terms containing a $d\bar z_j$ with 
$q+1\le j \le n$.  For $z''\in \Omega_2$ fixed, we can think of the
forms $\alpha_j$ and $\widehat g_j$ as $(0, q)$ and $(0, q-1)$-forms
respectively, in the variables $z'=(z_1,\ldots, z_q)$, $|z'|<1$.  Note
that we still have $\dbar_{z'} \widehat g_j =\alpha_j$,
where $\dbar_{z'}$ denotes $\dbar$ in the variables $z'$.
Let $\langle\cdot, \cdot\rangle$ be the standard pointwise inner
product on forms in $\C^q$. 
Let $\chi\in C^{\infty}_0(-\infty, \infty)$ be a cut-off function
such that $0\le \chi \le 1, \chi=1$ when $t\le 1/2$ and
$\chi =0$ when $t\ge 3/4$.  
Let $\beta=\chi(|z'|) d\bar z_1\wedge\cdots\wedge d\bar z_q$. 
It follows from the mean value property of holomorphic functions
that for $z''\in \Omega_2$,
$$
\align
|f_j(z'')-f_k(z'')| &= C\left|\int_{|z'|<1} \langle \alpha_j-\alpha_k, 
\beta \rangle dV(z')\right| \\
&= C\left|\int_{|z'|<1} \langle \widehat g_j- \widehat g_k, 
\vartheta\beta \rangle dV(z')\right| \quad (\text{$\vartheta$ is the formal
adjoint of $\dbar_{z'}$}) \\
&\le C \left\{\int_{|z'|<1} |\widehat g_j - \widehat g_k|^2\, dV(z')
\right\}^{1\over 2}
\endalign
$$
Therefore, after integrating in $z''$,  
$$
\|f_j-f_k\|_{\Omega_2}\le 
C\|\widehat g_j - \widehat g_k\|_{L^2_{(0, q-1)}(\Omega)}\le
C\|g_j - g_k\|_{L^2_{(0, q-1)}(\Omega)}
\quad \text{ as } \ j, k\to \infty.
$$
Since $\{f_j\}$ has no subsequence that converges in $L^2(\Omega_2)$,
$\{g_j\}$ has no subsequence that converges in 
$L^2_{(0, q-1)}(\Omega)$, contradicting the compactness of $S_q$. This completes the proof that (1) implies (2) in Theorem 1.1.
\enddemo

\bigskip

\subhead 5. Further Remarks \endsubhead
\smallskip

1) The arguments in Section 4 can be localized by using suitable
cut-off functions in $z''$ as well; it is enough to control the geometry locally. One then needs a lemma to the following effect:
Let $U_1$ and $U_2$ be neighborhoods of a boundary point $p_0$ of a bounded pseudoconvex domain $\Omega$, $U_1
\subset\subset U_2$.  Then $K_\Omega (w, w)$ and 
$\|K_\Omega(\cdot, w)\|_{\Omega\cap U_2}$ are comparable, uniformly
for $w\in \Omega\cap U_1$.  This can be shown by applying the reproducing
property of $K_{\Omega\cap U_2}$ to $K_\Omega(\cdot, w)$ and using
that $K_{\Omega\cap U_2}(w, w)$ and $K_\Omega(w, w)$ are comparable.
 Also, compactness of the $\dbar$-Neumann problem
is a local property: if every boundary point has the property that a
compactness estimate holds for forms supported near the point,
then the $\dbar$-Neumann problem is compact. This shows that Theorem~1.1
holds on domains that are {\it locally} convexifiable. 

2) It is noteworthy that in the proof that (1) implies (2)
in Theorem 1.1, we have only used that there is a compact solution
operator to $\dbar$ on the $(0, q)$-forms with holomorphic coefficients.

3) On smooth bounded {\it convex} domains there is a hierarchy of regularity
for the  $\dbar$-Neumann problem which can be described in terms of the 
contact with the boundary of affine complex varieties.  $N_q$ is subelliptic
if and only if the order of contact with the boundary of $q$-dimensional
affine complex varieties is bounded from above ([Ca3], [Mc], [Y]); $N_q$
is compact if and only if 
the boundary contains no $q$-dimensional affine varieties
(Theorem 1.1); finally, $N_q$ is globally regular regardless of whether
or not $b\Omega$ contains analytic varieties ([BS1]).

4) We have stated our results for $(0, q)$-forms, rather than $(p, q)$-forms,
as the index $p$ plays no r\^{o}le in solving $\dbar$.

5) To prove compactness of $N_q$, we have used (the analogue, for 
$(0, q)$-forms, of) Property~(P), see Proposition 3.1.  Our work
shows that for {\it convex} domains, this property is actually
{\it equivalent} to compactness of $N_q$.  On general pseudoconvex
domains, Property (P) still seems to be the  only systematic
way to  derive compactness of the $\dbar$-Neumann problem,
but (to quote from [BS3]) ``it is not yet understood how much
room there is between Property (P) and compactness''.

6) In the proof of Lemma 4.1, we have used the fact that on a bounded
convex domain $K_\Omega(z, z)$ $\to \infty$ as $z\to b\Omega$.  While
this is sufficient for the proof of Lemma 4.1, it is interesting
to note that for convex domains, there is the  (optimal) lower
estimate $K_\Omega(z, z)\ge C /(\text{dist}(z, b\Omega))^2$.  This
can be shown by the Ohsawa-Takegoshi
extension theorem (see [OT], [O2]) and the fact that 
the estimate is true for bounded convex domains in the plane.

\enddocument